\documentclass[11pt]{amsart}
\usepackage{amssymb,latexsym}

\theoremstyle{plain}
\newtheorem{theorem}{Theorem}
\newtheorem{corollary}{Corollary}
\newtheorem{lemma}{Lemma}
\newtheorem{proposition}{Proposition}
\theoremstyle{definition}

\theoremstyle{remark}

\newcommand{\mat}[4]{\begin{pmatrix}#1&#2\\ #3&#4\end{pmatrix}}
\newcommand{\matt}[2]{\begin{pmatrix}#1\\ #2\end{pmatrix}}

\begin{document}
\title[Theta Correspondence]
      {Theta Correspondence for $U(1,1)$ and $U(2)$}
\author{Jitka Stehnova}
\address{Department of Mathematics\\
         Oklahoma State University\\
         Stillwater, OK 74078}
\email{jstehnov@math.okstate.edu}
\date{September 8, 2008}
\begin{abstract}
  In this paper, we parametrize certain irreducible supercuspidal representations of $U(1,1)$ and $U(2)$ via explicit induction data. The parametrization depends on traceless elements of negative valuation in a quadratic extension of base field. We use the lattice model of the Weil representation to determine which traceless elements are involved in the theta correspondence for reductive dual pair $U(1,1)$ and $U(2)$. 
\end{abstract}

\maketitle

\section{Introduction}

Let $F$ be a p-adic field with odd residual characteristic. Let $E$ be a quadratic extension of $F$. Let $D$ be the quaternion division algebra over $F$ equipped with the nondegenerate symmetric bilinear form defined by the norm map $N_{D/F}$. Let $(V,(,)_{2})$ be a two-dimensional skew-Hermitian space over $E$ with a hyperbolic pair as basis and let $U(1,1)$ be the corresponding isometry group. Let $(W,(,)'_{2})$ be a two-dimensional Hermitian vector space over $E$ with $W=D$ and let $U(2)$ be the corresponding isometry group. Then $(U(1,1),U(2))$ forms a reductive dual pair in the symplectic group $Sp(\mathcal{W})$ where $\mathcal{W}=V\otimes_{E}W$ is a nondegenerate symplectic vector space equipped with the tensor product of the forms. 

Let $\chi$ be a nontrivial additive character of $F$. Let $\omega_{\chi}$ be the corresponding Weil representation associated to the metaplectic cover of $Sp(\mathcal{W})$. By restricting the Weil representation, we obtain a correspondence between some irreducible admissible representations of the metaplectic cover of $U(1,1)$ and some irreducible admissible representations of the metaplectic cover of $U(2)$. This correspondence is known as the \emph{theta correspondence} or \emph{Howe duality}. It is known that the correspondence is one-to-one by R.Howe ~\cite{Moeglin} and J.-P.Waldspurger ~\cite{Waldspurger}. In case of unitary groups, there are very few known examples of explicit theta correspondence, such as $(U_{1},U_{1})$ by Moen ~\cite{Moen} and recent results for small unitary groups by Misaghian ~\cite{Misaghian}, Stuffelbeam ~\cite{Stuffelbeam} and Pan ~\cite{Pan}. In this paper, we focus on two-dimensional unitary groups $U(1,1)$ and $U(2)$.

The fundamental problem is to parametrize irreducible admissible supercuspidal representations and use the lattice model of the Weil representation to construct explicit vectors that lead to a description of the local theta correspondence for the dual pair. Our approach follows the parametrization methods first used by Kutzko in ~\cite{Kutzko2} and ~\cite{Kutzko3} for $GL_{2}(F)$, then by Manderscheid in ~\cite{Mander1} for $SL_{2}(F)$. The advantage of this approach is that explicit construction of supercuspidals yields exact parametrizing data. 

The contents of this paper are as follows. In section $2$, we set up our notation and investigate groups $G=U(1,1)$ and $G'=U(2)$. In section $3$, we outline the induction method that we use to construct relevant representations of $U(1,1)$ and $U(2)$. In section $4$,  we construct some irreducible admissible supercuspidal representations of $U(1,1)$ with explicit induction data. We associate characters to tracelesss elements of $E$ and proceed with parametrization methods of Kutzko and Manderscheid. We repeat the process in section $5$ to construct some irreducible admissible supercuspidal representations of $U(2)$. In section $6$, we outline the lattice model of the Weil representation necessary in order to construct explicit vectors leading to the description of the theta correspondence. In section $7$, we determine which traceless elements are involved in the correspondence and show that they belong to corresponding cojugacy classes in $GL_{2}$ and the division algebra.

\section{Notation and Structure of $U(1,1)$ and $U(2)$}

Let $F$ be a nonarchimedean local field of residual characteristic $p$ with $p$ odd. Let $\mathcal{O}_{F}$ be the ring of integers of $F$, let $P_{F}$ be the maximal ideal in $\mathcal{O}_{F}$, let $\varpi=\varpi_{F}$ be uniformizer of $F$ and let $k_{F}$ be the residue field $\mathcal{O}_{F}/P_{F}$ with cardinality $q$. Let $v_{F}$ denote the valuation of $F$ and $\mathcal{O}^{\times}_{F}$ the group of units of $F$-space. Let $E$ be a quadratic extension of $F$. Let $\mathcal{O}_{E},P_{E},k_{E},\varpi_{E},v_{E}$ and $\mathcal{O}^{\times}_{E}$ play the corresponding roles with respect to $E$. 

Let $\tau: x\mapsto \bar{x}$ denote the nontrivial Galois automorphism of $E$ over $F$. Let $N=N_{E/F}$ and $Tr=Tr_{E/F}$ be the usual norm and trace maps associated to the extension $E/F$. Let $E^{1}$ be the group of norm-one elements in E and let $E^{0}$ denote the traceless elements of $E$. 

Let $D$ be the nonsplit quaternion division algebra over $F$ equipped with the nondegenerate symmetric bilinear form defined by the (reduced) norm map $N_{D/F}$. Let $\tau_{D}$ denote the involution on $D$ such that $N_{D/F}(x)=x\tau_{D}(x)$ and $Tr_{D/F}(x)=x+\tau_{D}(x)$. 

Let $a\in D, a\notin F$. The $F$-algebra $F[a]$ is then a field; since $D$ is an $F[a]$-vector space, we must have $[F[a]:F]=2$. Futher, there exists a separabale quadratic extension $E/F$ such that $E$ admits an $F$-embedding in $D$, see \cite{Bushnell}. More importantly, any quadratic field extension $E/F$ can be embedded in $D$. This has a profound effect. If $E$ is a quadratic subfield of $D$, then we may realize $D$ as the cyclic algebra $(E/F,\tau,a)$ where $a$ is an element of $F^{\times}$ which is not in the image of the norm map $N_{E/F}$. Notice that $\tau$ is precisely the restricion of the involution $\tau_{D}$ to $E$. In particular, given a generator $\alpha$ for $E/F$ there exists an element $\delta$ in $D^{\times}$ such that $\delta\alpha\delta^{-1}=\tau(\alpha)=-\alpha$ and $\delta^{2}=a$. We can take $\delta$ to be uniformizer of $D$ and let $\mathcal{O}_{D},P_{D},k_{D}$ play the corresponding roles. Then $\{1,\alpha,\delta,\alpha\delta\}$ forms a basis of $D$ over $F$ and $D=E\oplus\delta E$. 

If $E/F$ is unramified, we take $\varpi_{E}=\varpi=\delta^{2}$ and we choose an element $z$ of $F$, which is not a square in $F$, so that $E=F(\alpha)$ with $\alpha^{2}=z$,$\alpha\in\mathcal{O}^{\times}_{E}$ and $\alpha^{2}\in\mathcal{O}^{\times}_{F}$. We also have $|k_{E}|=q^{2}$. 

We realize $G=U(1,1)$ as the isometry group of a two-dimensional skew-Hermitian space $V$ over $E$ having a hyperbolic pair $\{\emph{u,v}\}$ as basis. In particular,
\newline
\begin{center}
$U(1,1)=\left\{x\in GL_{2}(E):\mat{0}{1}{-1}{0}=x\mat{0}{1}{-1}{0}\bar{x}^{t}\right\}$
\end{center}

\noindent For such $x$, $\det x.(\overline{\det x}) =1$. An easy argument reveals that $G=SU(1,1)\rtimes E^{1}$ where $SU(1,1)$ consists of the determinant-one elements of $G$. There is a natural identification of $SU(1,1)$ and $SL_{2}(F)$ and henceforth we consider $G=SL_{2}(F)\rtimes E^{1}$. 

We have three distinct copies of $E^{1}$ in G. The first is the semidirect copy, the second is the natural embedding of $E^{1}$ into $SL_{2}(F)$ (from a natural embedding of $E^{\times}$ into $GL_{2}(F)$) via 
\begin{center}
$x+y\alpha\mapsto\mat{x}{y}{y\alpha^{2}}{x}$.
\end{center}

\noindent The third copy of $E^{1}$ is the copy in the center of $G$ given by $\lambda\mapsto (\lambda,\lambda^{2})$. 

For $G'=U(2)$, we consider the Hermitian form on a two dimensional $E$-vector space $W=D$ defined by:

$(u,v)^{'}_{2}=\frac{1}{2}Tr_{D/E}(u\bar{v}),\  u=u_{1}+u_{2}\delta, \ v=v_{1}+v_{2}\delta , \ u_{1},u_{2},v_{1}, v_{2}\in E$.

\noindent It is easy to verify that the form is nondegenerate and anisotropic. We take $G'$ to be an isometry group of $(W,<,>')$.

\begin{center}
$U(2)=\left\{x\in GL_{2}(E):\mat{1}{0}{0}{-\delta^{2}}=x\mat{1}{0}{0}{-\delta^{2}}\bar{x}^{t}\right\}$.
\end{center}

\noindent Again, there is an important copy of $E^{1}$ in $G'$, namely $\mat{1}{0}{0}{\lambda}, \ \lambda\in E^{1}$. The anisotropic $U(2)$ has a subgroup of $SU(2)$ consisting of determinant-one elements and we can identify it with the group of norm-one elements in $D^{\times}$, $D^{1}=\{x=a+c\delta: a,c\in E, N_{D/F}(x)=1\}$. Similar argument as above reveals that $G'=D^{1}\rtimes E^{1}$.

\section{Irreducible Admissible Representations of $G,G'$}

To obtain irreducible admissible representations of $G$ and $G'$, we use the following approach. First, we associate characters of $F$ to traceless elements and extend them to characters of certain open compact subgroups of $G$ and $G'$. Second, we induce these characters using open compact induction to obtain irreducible admissible supercuspidal representations of $G,G'$. We make a frequent use of the following facts. By Jacquet ~\cite{Jacquet}, any irreducible smooth representation of a reductive $p$-adic group is admissible. And second, any irreducible representation obtained via compact induction from open compact subgroup is supercuspidal (~\cite{Bushnell2}). Thus, if we exhibit irreducibility of a compactly-induced representation, admissibility and supercuspidality automatically follow. The irreducible admissible representations of $G$ and $G'$ were previously studied by Stuffelbeam ~\cite{Stuffelbeam} and Misaghian ~\cite{Misaghian} and hence we omit some of the proofs. 

\section{Irreducible Admissible Representations of $U(1,1)$}

Recall $E^{1}$ is the group of are norm-one elements in $E^{\times}$ and $E^{0}$ is the group of the traceless elements in $E$. For any $r$, let $P^{r}_{E}=\{x\in E: x=a\varpi^{r}_{E},$ for some $a\in \mathcal{O}_{E}\}$ and define $P^{r}_{F}$ similarly. For $r\geq 1$, let $E^{1}_{r}=\{\lambda\in E^{1}:\lambda-1\in P^{r}_{E}\}$. We start with some important subgroups of $U(1,1)$. 

Henceforth, we assume $E/F$ is unramified. The subgroup $K_{1}=SL_{2}(\mathcal{O}_{F})\rtimes E^{1}$ is a maximal compact open subgroup of $U(1,1)$. Let
\begin{center}
$w=\mat{0}{1}{\varpi}{0}$
\end{center}

Then the other maximal compact open subgroup is, up to conjugacy, $K_{2}=w^{-1}K_{1}w$. Let $SL^{r}_{2}(\mathcal{O}_{F})=\{x\in SL_{2}(\mathcal{O}_{F}):x-1\equiv 0 \pmod{P^{r}_{F}}\}$. Let $K^{r}_{1}=E^{1}SL^{r}_{2}(\mathcal{O}_{F})\rtimes E^{1}$ and $K^{r}_{2}=w^{-1}K^{r}_{1}w$. These are again open compact subgroups of $U(1,1)$. 

Fix a nontrivial additive character $\psi$ of $F$ with the conductor $P_{F}$. Since $E=F(\alpha)$, we can naturally identify End$_{F}(E)$ with $M_{2}(F)$ via the map $a+b\alpha\mapsto\scriptsize{\mat{a}{b\alpha^{2}}{b}{a}}, a,b\in F$. Under this map, a traceless element is identified with a traceless matrix. Denote $M\_(F)$ the traceless matrices in $M_{2}(F)$. In particular for $\beta=y\alpha\in E^{0}$, we have a natural injection $E^{0}\rightarrow M\_(F)$ via
\begin{center}
$\beta=\mat{0}{y\alpha^{2}}{y}{0}$
\end{center} 

\noindent Let $M^{r}\_(F)$ denote the traceless matrices whose entries all have valuation $\geq r$. Let $n\geq 1$ and $r=\lfloor\frac{n+2}{2}\rfloor$. The map $x\mapsto x-1$ gives an isomophism from $M^{r}_{\_}(F)/M^{n+1}_{\_}(F)\cong SL^{r}_{2}(\mathcal{O}_F)/SL^{n+1}_{2}(\mathcal{O}_{F})$. This gives an isomorphism of character groups
\begin{center}
$M^{-n}_{\_}(F)/M^{1-r}_{\_}(F)\cong (SL^{r}_{2}(\mathcal{O}_{F})/SL^{n+1}_{2}(\mathcal{O}_{F}))\hat{}$
\end{center}
\noindent which yields a map $b\mapsto \psi_{b}$ where $\psi_{b}(x)=\psi(Tr(b(x-1)))$, see \cite{Bushnell}.

\begin{proposition}
Let $\beta\in E^{0}$ and $v_{E}(\beta)=-n$, $n\geq 1$. For $r=\lfloor\frac{n+2}{2}\rfloor$, we have $\psi_{\beta}:SL^{r}_{2}(\mathcal{O}_{F})\rightarrow\mathbb{C}^{\times}$ is a character trivial on $SL^{n+1}_{2}(\mathcal{O}_{F})$.
\end{proposition}
\begin{proof}
  This follows directly from the duality above.
\end{proof}

Next, we want $\psi_{\beta}$ is invariant under the semidirect action of $E^{1}$. From now on, fix $\beta\in E^{0}$ and $v_{E}(\beta)=-n$, $n\geq 1$.
\begin{proposition}
Let $x\in SL^{r}_{2}(\mathcal{O}_{F})$. We have $\psi_{\beta}(\sigma_{\lambda}(x))=\psi_{\beta}(x)$.
\end{proposition}
\begin{proof}
See ~\cite{Stehnova}, ~\cite{Stuffelbeam}.
\end{proof}

Let $SL_{2}(\mathcal{O}_{F})$ act on $\psi_{\beta}$ by conjugation: $\psi^{g}_{\beta}(h)=\psi_{\beta}(ghg^{1})$.

\begin{lemma}
The kernel of the action of $SL_{2}(\mathcal{O}_{F})$ on $\psi_{\beta}$ is $E^{1}SL^{n-r+1}_{2}(\mathcal{O}_{F})$.
\end{lemma}
\begin{proof}
See ~\cite{Stehnova}, ~\cite{Stuffelbeam}.
\end{proof}

For $\beta\in E^{0}$, let $\Lambda_{\beta}=\{\phi\in (E^{1})\hat{}:\phi\equiv \psi_{\beta}$ on $E^{1}\cap SL^{r}_{2}(\mathcal{O}_{F})\}$. If $\phi\in\Lambda_{\beta}$, then the map $\phi_{\beta}:E^{1}SL^{r}_{2}(\mathcal{O}_{F})\rightarrow \mathbb{C}^{\times}$ given by $\phi_{\beta}(\lambda x)=\phi(\lambda)\psi_{\beta}(x)$ is a well defined character. Let $\eta\in(E^{1})\hat{}$ and consider $\phi_{(\beta,\eta)}:K^{r}_{1}\rightarrow \mathbb{C}^{\times}$ given by $\phi_{(\beta,\eta)}(g,\gamma)=\phi_{\beta}(g)\eta(\gamma)$. This is again well defined character.

Again, we can conjugate these characters by $w$ and obtained well defined characters $\phi^{w}_{(\beta,\eta)}:K^{r}_{2}\rightarrow \mathbb{C}^{\times}$ where $\phi^{w}_{(\beta,\eta)}(\omega^{-1}(g,\gamma)\omega)=\phi_{(\beta,\eta)}(g,\gamma)$. 

In order to construct irreducible admissible supercuspidal representations of $U(1,1)$, we have to place further conditions on $v_{E}(\beta)=-n$ such as $n$ odd or even. If $n=2m+1$, then $r=m+1$ and $n-r+1=2m+1-m-1+1=m+1=r$ and hence the stabilizer of $\psi_{\beta}$ under the action of $SL_{2}(\mathcal{O}_{F})$ agrees with its domain on the $SL^{*}_{2}(\mathcal{O}_{F})$-part. If $n=2m$, then $r=m+1$ and $n-r+1=m=r-1$ so that the stabilizer of $\psi_{\beta}$ under the action of $SL_{2}(\mathcal{O}_{F})$ is strictly bigger on $SL^{*}_{2}(\mathcal{O}_{F})$-part than its group of definition.

\subsection{Odd Valuation}
In this section, we let $\beta\in E^{0}$ with $v_E(\beta)=-n=-(2m+1)$, $r=m+1$ and let $G=U(1,1)$. For $\eta\in(E^{1})\hat{}$, consider $\pi_{(\beta,\phi,\eta)}= $Ind$(G,K^{r}_{1};\phi_{(\beta,\eta)})$. 

\begin{theorem}
$\pi_{(\beta,\phi,\eta)}$ is an irreducible supercuspidal representation of $G$ with central character $\phi\cdot\eta^{2}$.
\end{theorem}
\begin{proof}
As mentioned above, we only need to show that $\pi_{(\beta,\phi,\eta)}$ is irreducible. Given a set of double coset representatives $\{x_{i}\}_{i\in I}$ for $E^{1}SL^{r}_{2}(\mathcal{O}_{F})\backslash SL_{2}(F)/E^{1}SL^{r}_{2}(\mathcal{O}_{F})$, it is clear that 
\begin{center}
$G=\bigcup_{i\in I}K^{r}_{1}\backslash (x_{i},1)/K^{r}_{1}$
\end{center}
\noindent Since the semidirect action of $E^{1}$ may have put some representatives in the same double coset, we will take an index subset $S\subset I$ such that $S$ has one representative for each double coset. Then by Mackey theory,
\begin{center}
$I(\pi_{(\beta,\phi,\eta)},\pi_{(\beta,\phi,\eta)})\cong\bigoplus_{i\in S}I(\phi_{(\beta,\eta)},\phi^{x_{i}}_{(\beta,\eta)})$
\end{center}

\noindent Upon the restriction to $E^{1}SL^{r}_{2}(\mathcal{O}_{F})\cap x^{-1}_{i}(E^{1}SL^{r}_{2}(\mathcal{O}_{F}))x_{i}$, $\phi_{(\beta,\eta)}=\phi_{\beta}$ and $\phi^{x_{i}}_{(\beta,\eta)}=\phi^{x_{i}}_{\beta}$. We also know by \cite{Mander3}, that the intertwining of $x_{i}$ outside the compact subgroup is $0$. Thus for each $i$,
\begin{center}
$I(\phi_{(\beta,\eta)},\phi^{x_{i}}_{(\beta,\eta)})\subset I(\phi_{\beta},\phi^{x_{i}}_{\beta})$
\end{center}

\noindent Since $n$ is odd, the stabilizer of $\psi_{\beta}$ is $E^{1}SL^{n-r+1}_{2}(\mathcal{O}_{F})=E^{1}SL^{r}_{2}(\mathcal{O}_{F})$. 
Hence by \cite{Mander1}, the representation Ind$(SL_{2}(F),E^{1}SL^{r}_{2}(\mathcal{O}_{F});\phi_{\beta})$ is an irreducible supercuspidal of $SL_{2}(F)$. Then Mackey theory and the selection of representatives give $I(\phi_{\beta},\phi^{x_{i}}_{\beta})=0$ for all $x_{i}\neq 1$. Since $I(\phi_{(\beta,\eta)},\phi_{(\beta,\eta)})=\mathbb{C}$, we conclude that $I(\pi_{(\beta,\phi,\eta)},\pi_{(\beta,\phi,\eta)})=\mathbb{C}$ and then \cite{Bushnell2}, Theorem 11.4  gives us $\pi_{(\beta,\phi,\eta)}$ is an irreducible representation. It follows it is admissible and supercuspidal. The statement about the central character follows from the definition of $\pi_{(\beta,\phi,\eta)}$.
\end{proof}

With the same notation as above, consider the character $\phi^{w}_{(\beta,\eta)}$ on $K^{r}_{2}$. Let $\pi'_{(\beta,\phi,\eta)}=$ Ind$(G,K^{r}_{2};\phi^{\omega}_{(\beta,\eta)})$.

\begin{theorem}
$\pi'_{(\beta,\phi,\eta)}$ is an irreducible supercuspidal representation of $G$ with the central character $\phi\cdot\eta^{2}$.
\end{theorem}
\begin{proof}
The construction is virtually the same, for details see ~\cite{Stehnova}, ~\cite{Stuffelbeam}.
\end{proof} 

\subsection{Heisenberg Construction}

In this section, we have $\beta\in E^{0}$ with $v_{E}(\beta)=-n=-2m, r=m+1$ and again $G=U(1,1)$. Since $n-r+1=m=r-1$, the character $\psi_{\beta}$ on $SL^{r}_{2}(\mathcal{O}_{F})$ is stabilized by subgroup $E^{1}SL^{r-1}_{2}(\mathcal{O}_{F})$. Hence for any $\eta\in(E^{1})\hat{}$, $I(G,K^{r}_{1};\phi_{(\beta,\eta)})$ will be reducible and we have to use different methods to find irreducible supercuspidal representations of $G$. The construction has been studied in ~\cite{Stuffelbeam}, therefore we only state the results. 

Let $SL^{r\_}_{2}(\mathcal{O}_{F})$ be the subset of $SL_{2}(\mathcal{O}_{F})$ such that the diagonal elements are congruent to $1$ modulo $P^{r-1}_{F}$ and off-diagonal elements are congruent to $0$ mod $P^{r}_{F}$. Then it is clear that $SL^{r}_{2}(\mathcal{O}_{F})\subset SL^{r\_}_{2}(\mathcal{O}_{F})\subset SL^{r-1}_{2}(\mathcal{O}_{F})$. Also, the character $\psi_{\beta}$ can be extended onto $SL^{r\_}_{2}(\mathcal{O}_{F})$ since $\beta\in E^{0}$ and hence $\psi_{\beta}$ depends only on off-diagonal elements.

Let $E^{1}_{0}=E^{1}\cap F^{\times}(1+P_{E})$. Then one checks that $E^{1}_{0}$  normalizes $SL^{r\_}_{2}(\mathcal{O}_{F})$. Depending on $\alpha$, $E^{1}$ may not normalize the above. Select $\phi\in(E^{1})\hat{}$ such that $\phi=\psi_{\beta}$ on $E^{1}\cap SL^{r\_}_{2}(\mathcal{O}_{F})$. Define $\phi_{\beta}: E^{1}_{0}SL^{r\_}_{2}(\mathcal{O}_{F})\rightarrow \mathbb{C}^{\times}$ naturally. It is clear that $\phi_{\beta}$ is a character. Now, when we add our semidirect product action, we may not necessarily obtain a group. Therefore instead of having $E^{1}$ in a semidirect product, we will work with $E^{1}_{1}=\{\lambda\in E^{1}:\lambda -1\in P^{1}_{E}\}$. This will guarantee us that $E^{1}_{0}SL^{r\_}_{2}(\mathcal{O}_{F})\rtimes E^{1}_{1}$ is a subgroup.

The extended $\psi_{\beta}$ is invariant under the semidirect action. And the kernel of the $SL_{2}(\mathcal{O}_{F})$-conjugate action on this extended $\psi_{\beta}$ is $E^{1}_{0}SL^{r_{\_}}_{2}(\mathcal{O}_{F})$. Let $\eta\in(E^{1})\hat{}$ and consider the character $\phi_{(\beta,\eta)}:E^{1}_{0}SL^{r\_}_{2}(\mathcal{O}_{F})\rtimes E^{1}_{1}\rightarrow \mathbb{C}^{\times}$ given by $\phi_{(\beta,\eta)}(\lambda g,\gamma)=\phi_{\beta}(\lambda g)\eta(\gamma)$ where $\eta$ is restricted to $E^{1}_{1}$. For later computations, we need exact number of matrices that form a complete set of distinct coset representatives for various cosets.

\begin{lemma}
$|E^{1}_{0}SL^{r-1}_{2}(\mathcal{O}_{F})\rtimes E^{1}_{1}:E^{1}_{0}SL^{r\_}_{2}(\mathcal{O}_{F})\rtimes E^{1}_{1}|=q$;

$|E^{1}_{0}SL^{r-1}_{2}(\mathcal{O}_{F})\rtimes E^{1}_{1}:E^{1}_{0}SL^{r}_{2}(\mathcal{O}_{F})\rtimes E^{1}_{1}|=q^{2}$;

$|E^{1}_{0}SL^{r}_{2}(\mathcal{O}_{F})\rtimes E^{1}_{1}\backslash E^{1}_{0}SL^{r-1}_{2}(\mathcal{O}_{F})\rtimes E^{1}_{1}/E^{1}_{0}SL^{r\_}_{2}(\mathcal{O}_{F})\rtimes E^{1}_{1}|=q$;

$|E^{1}_{0}SL^{r}_{2}(\mathcal{O}_{F})\rtimes E^{1}_{1}\backslash E^{1}SL^{r-1}_{2}(\mathcal{O}_{F})\rtimes E^{1}_{1}/E^{1}_{0}SL^{r}_{2}(\mathcal{O}_{F})\rtimes E^{1}_{1}|=\frac{q^{2}(q+1)}{2}$;

$|E^{1}SL^{r}_{2}(\mathcal{O}_{F})\rtimes E^{1}_{1}\backslash E^{1}SL^{r-1}_{2}(\mathcal{O}_{F})\rtimes E^{1}_{1}/E^{1}_{0}SL^{r}_{2}(\mathcal{O}_{F})\rtimes E^{1}_{1}|=q^{2}$;

$|E^{1}SL^{r}_{2}(\mathcal{O}_{F})\rtimes E^{1}_{1}\backslash E^{1}SL^{r-1}_{2}(\mathcal{O}_{F})\rtimes E^{1}_{1}/E^{1}SL^{r}_{2}(\mathcal{O}_{F})\rtimes E^{1}_{1}|=2q-1$
\end{lemma}
\begin{proof}
See ~\cite{Stehnova}.
\end{proof}

Consider $\rho^{o}_{(\beta,\phi,\eta)}=$Ind$(E^{1}_{0}SL^{r-1}_{2}(\mathcal{O}_{F})\rtimes E^{1}_{1}, E^{1}_{0}SL^{r\_}_{2}(\mathcal{O}_{F})\rtimes E^{1}_{1}; \phi_{(\beta,\eta)})$. This is a q-dimensional irreducible representation.

Consider $\rho_{(\beta,\phi,\eta)}=$Ind$(E^{1}_{0}SL^{r-1}_{2}(\mathcal{O}_{F})\rtimes E^{1}_{1}, E^{1}_{0}SL^{r}_{2}(\mathcal{O}_{F})\rtimes E^{1}_{1}; \phi_{(\beta,\eta)})$. This representation decomposes with the respect to $\rho^{o}_{(\beta,\phi,\eta)}$. By Lemma 2, we have exactly $q$ copies. 

Let $\tau_{(\beta,\phi,\eta)}=$Ind$(E^{1}SL^{r-1}_{2}(\mathcal{O}_{F})\rtimes E^{1}_{1}, E^{1}SL^{r}_{2}(\mathcal{O}_{F})\rtimes E^{1}_{1};\phi_{(\beta,\eta)})$.

\begin{proposition}
Let f be the character of $\tau_{(\beta,\phi,\eta)}$ and g the character of  Ind $(E^{1}SL^{r-1}_{2}(\mathcal{O}_{F})\rtimes E^{1}_{1}, E^{1}_{0}SL^{r}_{2}(\mathcal{O}_{F})\rtimes E^{1}_{1};\phi_{(\beta,\eta)})$. Then $2q^{-1}g-f$ is the character of a q-dimensional irreducible representation $\tau^{1}_{(\beta,\phi,\eta)}$ of $E^{1}SL^{r-1}_{2}(\mathcal{O}_{F})\rtimes E^{1}_{1}$ whose restriction to $E^{1}_{0}SL^{r-1}_{2}(\mathcal{O}_{F})\rtimes E^{1}_{1}$ is $\rho^{o}_{(\beta,\phi,\eta)}$.
\end{proposition}
\begin{proof}
For details, see ~\cite{Stehnova}, ~\cite{Stuffelbeam}.
\end{proof} 

\begin{lemma}
The representation $\tau^{1}_{(\beta,\phi,\eta)}$ extends to a unique q-dimensional irreducible representation of $K^{r-1}_{1}$.
\end{lemma}
\begin{proof}
For detailed construction, see ~\cite{Stuffelbeam}.
\end{proof}

Define $\pi_{(\beta,\phi,\eta)}=$Ind$(G,K^{r-1}_{1};\tau^{1}_{(\beta,\phi,\eta)})$.
\begin{theorem}
$\pi_{(\beta,\phi,\eta)}$ is an irreducible supercuspidal representation of $G$ with the central character $\phi\cdot\eta^{2}$.
\end{theorem}
\begin{proof}
The argument is analogous to that involved in proving Theorem 4.1.1.
\end{proof}

Keeping the same $\beta,\eta,r,\phi,n$, we will construct irreducible supercuspidals on $K^{r}_{2}$. Consider $\phi^{w}_{(\beta,\eta)}$ on $K^{r}_{2}$. We can use analogous arguments, properly modify them  and reproduce the unique q-dimensional irreducible representation $\tau^{1}$ of $K^{r-1}_{2}$. Hence we only list the following result. Let $\pi'_{(\beta,\phi,\eta)}=$Ind$(G,K^{r-1}_{2};\tau^{1}_{(\beta,\phi,\eta)})$. 
\begin{theorem}
$\pi'_{(\beta,\phi,\eta)}$ is an irreducible supercuspidal representation of $G$ with the central character $\phi\cdot\eta^{2}$.
\end{theorem}

\subsection{The Level Zero Case}

In ~\cite{Gerardin}, Gerardin defines the Weil representation for symplectic groups, general linear groups and unitary groups over finite fields. He canonically identifies $U(2,k_{E})\cong SL_{2}(k_{F})\rtimes k^{1}_{E}$ with a subgroup of $Sp(4,k_{F})$. The results are applicable only if $E/F$ is unramified. 

Let $\chi$ be an additive character of $k_{F}$ and $\omega_{\chi}$ be the associated Weil representation of $Sp(4,k_{F})$. Gerardin proves that the Weil representation restricted to $U(2,k_{E})$ decomposes into irreducibles,
\begin{center}
$\omega_{\chi|U(2,k_{E})}= sgn \otimes \left(\displaystyle{\bigoplus_{\xi\in(E^{1}/E^{1}_{1})\hat{}}} \vartheta_{\xi}\right)$
\end{center}

\noindent where sgn is the unique nontrivial quadratic character of $U(2,k_{E})$ and $\vartheta_{1}$ is a q-dimensional irreducible representation withe the central character $1$, and for $\xi\neq 1$, $\vartheta_{\xi}$ is a $(q-1)$-dimensional irreducible cuspidal representatin of $U(2,k_{E})$ with the central character $\xi$. 

We assume $\xi\neq 1\in (E^{1}/E^{1}_{1})\hat{}$ and use the corresponding cuspidal representation $\vartheta_{\xi}$ to construct irreducible admissible supercuspidals of $U(1,1)$. We may lift $\vartheta_{\xi}$ to an irreducible $(q-1)$-dimensional representation $\rho_{(\xi,\eta)}$ of $K_{1}$. The construction and results  are known, hence we only state the important theorems.

Define $\pi_{(\xi,\eta)}=$Ind$(G,K_{1};\rho_{(\xi,\eta)})$. 
\begin{theorem}
$\pi_{(\xi,\eta)}$ is an irreducible admissible supercuspidal representation of $G$ with central character $\xi\cdot\eta^{2}$.
\end{theorem}
\begin{proof}
See ~\cite{Stehnova, Stuffelbeam}.
\end{proof}

In similar manner, we construct supercuspidals from $K_{2}$. For a nontrivial character $\xi\in(E^{1}/E^{1}_{1})\hat{}$ and $\eta\in (E^{1})\hat{}$, the representation $\rho^{w}_{(\xi,\eta)}$ is irreducible on $K_{2}$. Define $\pi'_{(\xi,\eta)}=$Ind$(G,K_{2};\rho^{w}_{(\xi,\eta)})$. The induced representation is an irreducible supercuspidal of $G$ with central character $\xi\cdot\eta^{2}$. 

\section{Irreducible Admissible Representations of $U(2)$}

Recall $G'=U(2)=D^{1}\rtimes E^{1}$ where $D^{1}$ are norm-one elements of $D$. Also, recall that $\delta$ is uniformizer of $D$ with $\delta^{2}=\varpi$. For any $r$, let $P^{r}_{D}=\{x\in D: x=a\delta^{r}, a\in\mathcal{O}_{D}\}$. Let $D^{0}$ denote the traceless elements in $D$ and let $D^{1}_{r}=\{x\in D^{1}: x-1\in P^{r}_{D}\}$. Also notice that, due to ramification of $D$ over $F$, we have $F\cap P^{n}_{D}\subset P^{\lfloor\frac{n+1}{2}\rfloor}_{F}$. Since $G'$ is compact, all of its irreducible admissible representations are supercuspidal. The supercuspidal representation of $U(2)$ were previously studied in \cite{Misaghian}. Our approach provides the explicit construction of supercuspidals with exact parametrizing date. The construction is virtually the same as for case $U(1,1)$ and hence we only state the important results. 

\subsection{Characters of $D^{1}$} 
\noindent First, we look at one dimensional irreducible representations. Let $D^{1}_{1}=\{x\in D^{1}:x-1\in P_{D}\}$. It is a well known fact that the commutator group of $G'$, $[G',G']=D^{1}_{1}\rtimes \{1_{E}\}$. Straightforward computations show that $D^{1}/D^{1}_{1}$ is a cyclic group of order $q+1$. 
\begin{lemma}
There is a bijection between characters of $D^{1}$ and characters of $D^{1}/D^{1}_{1}$.
\end{lemma}
\begin{proof}
Clear.
\end{proof}

\subsection{Characters Associated to Traceless Elements}
This construction is virtually the same as in the case of $U(1,1)$. Fix a nontrivial additive character $\psi$ of F with the conductor $P_{F}$. Recall, $E/F$ is unramified.

\begin{proposition}
Let $\beta\in D^{0}$ and $v_{D}(\beta)=-n$, $n\geq 1$. For $r=\lfloor\frac{n+2}{2}\rfloor$, define $\psi_{\beta}: D^{1}_{r}\rightarrow \mathbb{C}^{\times}$ by $\psi_{\beta}(h)=\psi(Tr(\beta(h-1)))$, $h\in D^{1}_{r}$. Then $\psi_{\beta}$ is a character of $D^{1}_{r}$ trivial on $D^{1}_{n+1}$. 
\end{proposition}

Next we want to show $\psi_{\beta}$ is invariant under the semidirect action of $E^{1}$. From now on fix $\beta\in D^{0}$ and $v_{D}(\beta)=-n, n\geq 1$.
\begin{proposition}
Let $h\in D^{1}_{r}$. We have $\psi_{\beta}(\sigma'_{\lambda}(h))=\psi_{\beta}(h)$, for $\lambda\in E^{1}$.
\end{proposition}
\begin{proof}
Modify the proof of Proposition $2$, Section $4$. For details, see ~\cite{Stehnova}.
\end{proof}

For $g\in D^{1}$, define $\psi^{g}_{\beta}(h)=\psi_{\beta}(ghg^{-1})$. The action is well defined, since for chosen $g$, $ghg^{-1}\in D^{1}_{r}$ for $x\in D^{1}_{r}$. Thus we can determine the stabilizer.
\begin{lemma}
The stabilizer of the action of $D^{1}$ on $\psi_{\beta}$ is $E^{1}D^{1}_{n-r+1}$. 
\end{lemma}
\begin{proof}
Simiarly to the proof of Lemma $1$, Section $4$.
\end{proof}

For $\beta\in D^{0}$, let $\Lambda_{\beta}=\{\gamma\in (E^{1})\hat{}:\gamma\equiv\psi_{\beta}$ on $E^{1}\cap D^{1}_{r}\}$. If $\gamma\in\Lambda_{\beta}$, then the map $\gamma_{\beta}:E^{1}D^{1}_{r}\rightarrow\mathbb{C}^{\times}$ given by $\gamma_{\beta}(\lambda x)=\gamma(\lambda)\psi_{\beta}(x)$ is a well defined character. Let $\zeta\in (E^{1})\hat{}$ and consider $\gamma_{(\beta,\zeta)}:E^{1}D^{1}_{r}\rtimes E^{1}\rightarrow\mathbb{C}^{\times}$ given by $\gamma_{(\beta,\zeta)}(g,\lambda')=\gamma_{\beta}(g)\zeta(\lambda')$. This is again well defined character. 

In order to construct the irreducible admissible supercuspidal representations of $U(2)$, we have to place further conditions on $v_{D}(\beta)=-n$, such as $n$ odd or even. If $n=2m+1$, then $r=m+1$ and $n-r+1=2m+1-m-1+1=m+1=r$ and hence the stabilizer of $\varphi_{\beta}$ under the action of $D^{1}$  and its domain agrees on $D^{1}_{*}$-part. If $n=2m$, then $r=m+1$ and $n-r+1=m=r-1$ so that the stabilizer of $\psi_{\beta}$ under the action of $D^{1}$ could be bigger on the $D^{1}_{*}$-part than its group of definition. Thus we will have to place additional conditions on $r$, such as $r$ is odd or even.

\subsection{Odd Valuation}
Let $\beta\in D^{0}$ with $v_{D}(\beta)=-(2m+1)$, $r=m+1$, and let $G'=U(2)$. For $\zeta\in (E^{1})\hat{}$, consider $\pi_{(\beta,\gamma,\zeta)}=$Ind$(G',E^{1}D^{1}_{r}\rtimes E^{1};\gamma_{(\beta,\zeta)})$. 
\begin{theorem}
$\pi_{(\beta,\gamma,\zeta)}$ is an irreducible admissible representation of $G'$.
\end{theorem}
\begin{proof}
Given the set of double coset representatives $\{x_{i}\}_{i\in I}$ for $E^{1}D^{1}_{r}\backslash D^{1}/E^{1}D^{1}_{r}$, it is clear that
\begin{center}
$G'=\bigcup_{i\in I} E^{1}D^{1}_{r}\rtimes E^{1}\backslash (x_{i},1)/E^{1}D^{1}_{r}\rtimes E^{1}$
\end{center}

\noindent Since the semidirect action of $E^{1}$ may have related some representatives, we will take an index subset $J\subset I$ such that $J$ has one representative for each double coset. Then by Mackey Theory,
\begin{center}
$I(\pi_{(\beta,\gamma,\zeta)},\pi_{(\beta,\gamma,\zeta)})\cong\displaystyle{\bigoplus_{i\in J}}I(\gamma_{(\beta,\zeta)},\gamma^{x_{i}}_{(\beta,\zeta)})$
\end{center}

Upon the restriction to $E^{1}D^{1}_{r}\cap x^{-1}_{i}(E^{1}D^{1}_{r})x_{i}, \gamma_{(\beta,\zeta)}=\gamma_{\beta}$ and $\gamma^{x_{i}}_{(\beta,\zeta)}=\gamma^{x_{i}}_{\beta}$. Thus for each $i$,
\begin{center}
$I(\gamma_{(\beta,\zeta)},\gamma^{x_{i}}_{(\beta,\zeta)})\subset I(\gamma_{\beta},\gamma^{x_{i}}_{\beta})$
\end{center}

Since $n$ is odd, the stabilizer of $\psi_{\beta}$ is exactly the domain of it and hence by Clifford Theory, theorem (45.2)' in \cite{Curtis}, the representation Ind$(D^{1},E^{1}D^{1}_{r};\gamma_{\beta})$ is an irreducible supercuspidal of $D^{1}$. Then Mackey theory and the selection of representatives give $I(\gamma_{\beta},\gamma^{x_{i}}_{\beta})=0$ for $x_{i}\neq 1$. Since $I(\gamma_{(\beta,\zeta)},\gamma_{(\beta,\zeta)})=\mathbb{C}$, we conclude that $I(\pi_{(\beta,\gamma,\zeta)},\pi_{(\beta,\gamma,\zeta)})=\mathbb{C}$. By \cite{Bushnell}, Theorem 11.4 it follows that $\pi_{(\beta,\gamma,\zeta)}$ is irreducible admissible representation and since $U(2)$ compact, it is supercuspidal. 
\end{proof}

\subsection{Even Valuation}
In this section, we take $v_{D}(\beta)=-n=-(2m)$, then $r=m+1$ and $n-r+1=m=r-1$ so that the stabilizer of $\psi_{\beta}$ under the action of $D^{1}$ could be bigger than its domain. In order to find irreducible admissible representations, we have to place additional conditions on $r$, such as $r$ is odd or even. For this section, we assume $r$ is odd. It follows $m$ is even and $n$ is divisible by $4$. Let $D^{*}$ denote the set of elements $x$ in $D$ such that $1+x$ is invertible. Then recall that the Cayley transform is the well defined map $c$ from $D^{*}$ to itself defined by $c(x)=(1-x)(1+x)^{-1}$. Notice that $c$ is a bijection onto $D^{*}$ with inverse $c$ itself.  

\begin{lemma}
If $r$ odd, then 
\begin{center}
$(E^{1}D^{1}_{r-1})/D^{1}_{n+1}=(E^{1}D^{1}_{r})/D^{1}_{n+1}$
\end{center}
where $r-1=m=n/2$ and $r=m+1=n/2+1$.
\end{lemma}

\begin{proof}
Since $E/F$ is unramified, ${\bf{k}}_{E}={\bf{k}}_{F}$. Let $\scriptsize{h=\frac{(1-a\delta^{r-1})}{(1+a\delta^{r-1})}}D^{1}_{n+1}$ be an element of $D^{1}_{r-1}/D^{1}_{n+1}$. Notice that we can write $h$ in this form due to the Cayley transform. Write $a=a_{0}+a_{1}\delta$ where $a_{0},a_{1}\in\mathcal{O}_{E}$. Now we have:
\begin{align*}
h&=\frac{1-a\delta^{r-1}}{1+a\delta^{r-1}}D^{1}_{n+1}=\frac{1-(a_{0}+a_{1}\delta)\delta^{r-1}}{1+(a_{0}+a_{1}\delta)\delta^{r-1}}D^{1}_{n+1}\\
&=\left (\frac{1-a_{0}\delta^{r-1}}{1+a_{0}\delta^{r-1}}\right )\left (\frac{1-a_{1}\delta^{r}}{1+a_{1}\delta^{r}}\right )\left (\frac{1-(1-a_{0}\delta^{r-1})^{-1}(1-a_{1}\delta^{r})^{-1}(a_{1}a_{0}\delta^{n+1})}{1-(1+a_{0}\delta^{r-1})^{-1}(1+a_{1}\delta^{r})^{-1}(a_{1}a_{0}\delta^{n+1})}\right )
\end{align*}

\vspace{.5cm}
\noindent Since $r$ odd, $r-1=m$ even and hence $a_{0}\delta^{r-1}\in E$ and thus the first quotient is in $E^{1}$. By definition, the second quotient is in $D^{1}_{r}=D^{1}_{n/2+1}$. To obtain the result, it suffices to show that third quotient is in $D^{1}_{n+1}$. Quick calculations shows the quotient is in $D^{1}$ and after subtracting $1$, the quotient is $\equiv 0 \pmod {P^{n+1}_{D}}$. Thus 
\begin{center}
$(E^{1}D^{1}_{r-1})/D^{1}_{n+1}\subset (E^{1}D^{1}_{r})/D^{1}_{n+1}$
\end{center}

\noindent The other containment is true by definition of the filtration, and we are done.  
\end{proof}

Thus any character of $(E^{1}D^{1}_{r-1})/D^{1}_{n+1}$ is a character of $(E^{1}D^{1}_{r})/D^{1}_{n+1}$ and vice versa. Hence we can take $\psi_{\beta}$ as in Proposition $1$ , and obtain a character $\gamma_{(\beta,\zeta)}$. By above Lemma, they both are characters on $E^{1}D^{1}_{r-1}\rtimes E^{1}$ and hence the stabilizer of $\psi_{\beta}$ under the action of $D^{1}$ coincides with its domain. Consider $\pi'_{(\beta,\gamma,\zeta)}=$Ind$(G',E^{1}D^{1}_{r-1}\rtimes E^{1}; \gamma_{(\beta,\zeta)})$. 
\begin{theorem}
$\pi'_{(\beta,\gamma,\zeta)}$ is an irreducible admissible representation of $G'$.
\end{theorem}
\begin{proof}
Apply the proof of Theorem $4.3.1$ with $r=r-1$. 
\end{proof}  

\subsection{Even Valuation - Heisenberg Construction}

In this case, we have $n$-even, $r$-even. Hence the stabilizer of $\psi_{\beta}$ under the action of $D^{1}$ is $E^{1}D^{1}_{r-1}$, strictly bigger than its domain. The Heisenberg construction is the same as in case $U(1,1)$, hence we only state important results.  

\begin{lemma}
There is a unique $q$-dimensional irreducible representation $\tau^{1}_{(\beta,\gamma,\zeta)}$  of $E^{1}D^{1}_{r-1}\rtimes E^{1}$.
\end{lemma}
\begin{proof}
For details, see ~\cite{Stehnova}.
\end{proof}

Define $\pi'_{(\beta,\gamma,\zeta)}=$Ind$(G',E^{1}D^{1}_{r-1}\rtimes E^{1};\tau^{1}_{(\beta,\gamma,\zeta)})$.
\begin{theorem}
\label{irreducible2:evenh:thm1}
$\pi'_{(\beta,\gamma,\zeta)}$ is an irreducible representation of $G'$.
\end{theorem}
\begin{proof}
Simiarly to the proof of Theorem $3$, Section $4.2$.
\end{proof}

\section{Lattice Model of the Weil Representation}
In this section, we detail the method of lattice models of the Weil representation. This will allow is to explicitly determine the occurence of irreducible admissible representations of $U(1,1)$ and $U(2)$ as quotients of the smooth Weil representation. This section is a recapitulation of a material in \cite{Pan} and 
\cite{Mander1} modified to fit our needs in the next section.

Let $(\mathcal{W}, <,>)$ be a nondegenerate symplectic vector space of dimension $2n$ over $F$ and let $H(\mathcal{W})$ be the associtated Heisenberg group, $H(\mathcal{W})=\mathcal{W}\oplus F$ with $F$ being the center of $H(\mathcal{W})$. Let $\chi$ be a notrivial additive character of $F$ and let $\rho_{\chi}$ the associated unique unitary representation of $H(\mathcal{W})$ with central character $\chi$. Let $\omega_{\chi}$ be the correspoding Weil representation and denote $\omega^{\infty}_{\chi}, \rho^{\infty}_{\chi}$ corresponding smooth representations.

In this section, we obtain a realization of $\omega_{\chi}$ by working with certain non-self-dual lattices in $\mathcal{W}$. We begin recalling some features of the Weil representation over the finite field $k_{F}$. 

Let $\mathcal{W'}$ be a finite dimensional vector space over $k$. Suppose that $\mathcal{W'}$ is equipped with a nondegenerate skew-symmetric bilinear form $<,>$ and let $G(\mathcal{W'})$ be its isometric group. Let $H(\mathcal{W'})$ denote the Heisenberg group attached to $\mathcal{W'}$ and let $\chi$ be a nontrivial additive character of $k$. Then, there is a unique (up to equivalence) unitary representation of $\rho_{\chi}$ of $H(\mathcal{W'})$ with central character $\chi$. Moreover, there is a representation $\omega_{\chi}$ of $G$ on the space $\rho_{\chi}$ such that
\begin{center}
$\omega_{\chi}(g)\rho_{\chi}(h)=\rho_{\chi}(gh)\omega_{\chi}(g)$
\end{center}
\noindent where $h\in H(\mathcal{W'}),g\in G$. This representation is unique up to equivalence except in the case where $\text{dim}_{F}\mathcal{W'}=2$ and $|k|=3$. In this exceptional case, we may fix $\omega_{\chi}$ to satisfy a certain condition in a Schrodinger model. For more details, see \cite{Mander1}.

We now turn to the lattice model for a certain type of lattice that is not self-dual. We return to the notation of the previous section and suppose $L$ is an $\mathcal{O}_{F}$-lattice in $\mathcal{W'}$ which altough is not self-dual does satisfy 
\begin{center}
$P_{F}L^{*}\subseteq L\subsetneq L^{*}$
\end{center}

\noindent We will call this lattice a \emph{good lattice}. Notice that $\bar{L}=L^{*}/L$ is an even dimensional vector space over $k$. 

Let $d$ be an interger such that $\chi$ is trivial on $P^{d}_{F}$ but not trivial on $P^{d-1}_{F}$. Let $\bar{x}$ and $\bar{y}$ in $\bar{L}$ be a preimages of $x$ and $y$ in $L^{*}$ and set $<\bar{x},\bar{y}>_{d}=\varpi^{1-d}<x,y>$. One can check that $<,>_{d}$ is well-defined nondegenerate skew-symmetric bilinear form on $\bar{L}$. We may also define a character $\chi'$ of $k$ by setting $\chi'(\bar{x})=\chi(x)$ where $\bar{x}\in k$ and $x$ is an element of $P^{d-1}_{F}/P^{d}_{F}$ with image $\bar{x}$ under the map induced by $y\mapsto \varpi^{1-d}y$ from $P^{d-1}_{F}$ to $\mathcal{O}_{F}$. Let $\rho_{\chi'}$ denote a representation of $H(\bar{L})$ with central character $\chi'$. 

Let $J^{*}$ be the subgroup of $H(\mathcal{W'})$ generated by $e(L^{*})$ and let $J$ be the subgroup of $H(\mathcal{W'})$ generated by $e(L)$. Then we may inflate $\rho_{\chi'}$ to a representation of $J^{*}$ which is trivial on $J$. We also define $\rho_{L}$ a representation of $\gamma^{-1}(L^{*})$ on the space of $\rho_{\chi'}$ by $\rho_{L}(ah)v=\chi(a)\rho_{\chi'}(h)v$ where $a\in Z(H(\mathcal{W'})), h\in J^{*}, v$ in the space of $\rho_{\chi'}$. Then Ind$(H(\mathcal{W'}),\gamma^{-1}(L^{*});\rho_{L})$ realizes $\rho_{\chi}$. For more details and proofs, see \cite{Mander1}.

We now need to make this realization more explicit. Let $X$ be the finite dimensional Hilbert space of $\rho_{L}$ and let $||\ ||$ denote the norm on $X$. Let $Y$ denote the space of $\rho_{\chi}$ and let $S_{L}$ denote the set of coset representatives for $\mathcal{W'}/L^{*}$. Then $Y$ is the set of functions $f:\mathcal{W'}\rightarrow X$ satisfying:
\begin{enumerate}
  \item [(i)] $f(w+a)=\chi(<w,a>/2)\rho_{L}(e(a))f(w)$ for $a\in L^{*}$
  \item [(ii)] $\sum_{w\in S_{L}} ||f(w)||^{2}<\infty$
\end{enumerate}

\noindent The action of $\rho_{\chi}$ is given by 
\begin{center}
$(\rho_{\chi}(e(w))f)(w')=\chi(<w',w>/2)f(w'+w)$
\end{center}
\noindent for $f\in Y, w,w'\in\mathcal{W'}$. For each $w\in\mathcal{W'}, x\in X$ of length one, let $y_{w,x}$ denote the function on $Y$ supported on $-w+L^{*}$ taking the value $x$ at $-w$. Then, if we choose an orthonormal basis $S_{X}$ for $X$, we have that $Y$ consists of linear combinations 
\begin{center}
$\displaystyle{\sum_{\begin{subarray}{1}
        w\in S_{L}\\
        w\in S_{X}
       \end{subarray}} }
  a_{w,x}y_{w,x}$
\end{center}
\noindent with
\begin{center}
$\sum |a_{w,x}|^{2}<\infty $
\end{center}
\noindent and $Y^{\infty}$ is the subspace of $Y$ consisting of finite linear combinations of the above form. 

We now consider  $\omega_{\chi}$. Let $K$ be the maximal compact subgroup of $Sp(\mathcal{W'})$ which stabilizes $L^{*}$ and let $K'$ be a subgroup of $K$ acting trivially on $L^{*}/L$. We may identify $K/K'$ with the isometry group of the symplectic space $\bar{L}$ and thus there exists a unique representation $\omega_{\chi}$ of $K$ on $X$ which is trivial on $K'$ and satisfies
\begin{center}
$\omega_{L}(g)\rho_{L}(h)=\rho_{L}(gh)\omega_{L}(g)$
\end{center}

\noindent for $h$ in $\gamma^{-1}(L^{*})$ and $g\in K$. 

\begin{proposition}
The representation $\omega_{\chi}$ may be chosen so that it restricts to a representation of $K$. In particular, $\omega_{\chi}$ may be chosen so that for $f$ in $Y$ and $k$ in $K$
\begin{center}
$\omega_{\chi}(k)f(w)=\omega_{L}(k)f(k^{-1}w)$
\end{center}
\noindent and thus for $k$ in $K$ 
\begin{center}
$\omega_{\chi}(k)y_{w,x}=y_{kw,\omega_{L}(k)x}$.
\end{center}
\noindent In addition, the space of smooth vectors $Y^{\infty}$ for $\omega_{\chi}$ consists of those $f$ in $Y$ supported on a finite number of $\ \mathcal{W'}/L^{*}$ cosets, i.e. those $f$ which are finite linear combinations of the $\{y_{w,x}\}$.
\end{proposition}
\begin{proof}
See, for example, Chapter 5 of \cite{Moeglin}.
\end{proof}

Now suppose, $L$ is a lattice in $\mathcal{W'}$ as above and $M$ is a sublattice of $L$. Then $H_{M}=\{g\in G|(g-1)M^{*}\subset L^{*}\}$ is a subgroup of $G$. And futher, we have

\begin{proposition}
If a function $f$ in $Y$ is supported on $M^{*}$, then
\begin{center}
$\omega_{\chi}(h)f(w)=\rho_{L}(2c(h)w)\chi(<w,c(h)w>)f(w)$
\end{center}
\noindent for $h$ in $H_{M}$ where $c(h)=(1-h)(1+h)^{-1}$ is the Cayley transform of $h$.
\end{proposition}
\begin{proof}
This result can be proved with a straightforward modification of the proof of \cite{Moeglin}.
\end{proof}

\section{Theta Correspondence}
In this section, we use the methods of previous section to begin to determine which irreducible admissible representations of $U(1,1)$ and $U(2)$ occur in the theta correspondence. For additive character of $\psi$ fixed in section $4,5$, set $\chi=\psi_{\varpi}$, that is, $\chi(x)=\psi_{\varpi}(x)=\psi(\varpi x)$ for $x$ in $F$. Then $\psi\cdot Tr_{E/F}$ is a character of $E$ with conductor $P_{E}$. Recall $\mathcal{W}=V\otimes_{E} W$ is equipped with a nondegenerate skew-symmetric bilinear form $<<,>>$ by setting $<<v_{1}\otimes w_{1},v_{2}\otimes w_{2}>>=Tr(<v_{1},v_{2}>\overline{<w_{1},w_{2}>'})$. 

Also recall $G=U(1,1)$ is the isometry group of $<,>_{2}$, $G'=U(2)$ is the isometry group of $<,>'_{2}$. We may identify $G$ and $G'$ with subgroups of $Sp(8)$ by letting $G$ act on $\mathcal{W}$ by premultiplication by inverses and letting $G'$ act on $\mathcal{W}$ by postmultiplication. Note that in this identification $G$ and $G'$ are each other commutants in $Sp(8)$, i.e. form a reductive dual pair. Recall $E/F$ is unramified. 

Let $\Gamma=\mathcal{O}_{E}u+\mathcal{O}_{E}v$ be the lattice in $V$ and $\Gamma'=\mathcal{O}_{E}+\mathcal{O}_{E}\delta$ be the lattice in $W$. Then $A= \Gamma\otimes \Gamma'$ is a lattice in $\mathcal{W}$. 

\begin{lemma}
$A$ is a non-self dual lattice and $A^{*}=(\mathcal{O}_{E}u+\mathcal{O}_{E}v)\otimes (\mathcal{O}_{E}+ P^{-1}_{E}\delta)$.
\end{lemma}
\begin{proof}
This can be checked directly, see for example ~\cite{Pan}, ~\cite{Stehnova}.
\end{proof}

Notice that $A$ is a non-self dual ``good'' lattice, i.e. satisfying 
\begin{center}
$\varpi_{F} A^{*}\subseteq A\subset A^{*}$
\end{center}
\noindent and hence we may apply results detailed in previous section. Let $A_{F}(\mathcal{W})=Hom_{F}(\mathcal{W},\mathcal{W})$ and for $k\in\mathbb{Z}, M^{k}=P^{k}A$. Let $M=\{M^{k}\}_{k\in\mathbb{Z}}$ be a lattice chain and let $\mathcal{A}$ be the subring of $A_{F}(\mathcal{W})$ consisting of elements $x$ such that $xM^{k}\subseteq M^{k}$ for all $k$. Also, for $n\geq 1$, let $\mathcal{P}^{n}$ be the set of elements $x$ in $\mathcal{A}$ satisfying $xM^{k}\subseteq M^{k+n}$ for all $k$. Let $U(\mathcal{A})=\{x\in Sp(\mathcal{W}):x\in\mathcal{A}^{\times}\}$ and for $n\geq 1, U^{n}(\mathcal{A})=\{x\in Sp(\mathcal{W}): x-1\in\mathcal{P}^{n}\}$. Finally, let $U^{n}_{1}(\mathcal{A})=U^{n}(\mathcal{A})\cap G$ and $U^{n}_{2}(\mathcal{A})=U^{n}(\mathcal{A})\cap G'$. Notice that these filtrations correspond to the filtrations on $G$ and $G'$ defined in Sections $4$ and $5$. 

\begin{lemma}
For $k\geq 0$, $(M^{k})^{*}=P^{-k}A^{*}$.
\end{lemma}
\begin{proof}
Recall $M^{k}=P^{k}A$, hence $M^{k}=(P^{k}\Gamma)\otimes\Gamma'$ or $M^{k}=\Gamma\otimes(P^{k}\Gamma')$. In the first case one can check, $(M^{k})^{*}=(P^{-k}\Gamma)\otimes (\Gamma')^{*}=P^{-k}A^{*}$. In the latter case, $(M^{k})^{*}=\Gamma\otimes(P^{-k})(\Gamma')^{*}=P^{-k}(\Gamma\otimes(\Gamma')^{*})=P^{-k}A^{*}$.
\end{proof}

Now to use a result of Section 5.2, we fix a set $S_{A}$ of coset representatives for $\mathcal{W}/A^{*}$. Recall for $v\in V, w\in W, x\in X$, $y_{v\otimes w}, x$ denotes the function $f$ in $Y$ supported on $-(v\otimes w)+A^{*}$ and taking the value $x$ at $-(v\otimes w)$. Also recall $\bar{A}=A^{*}/A$ is a 4-dimensional vector space over $k_{F}$. 

\begin{lemma}
\label{theta:unram:lemma3}
Let $v,v'\in V, w,w'\in W$ and $x,x'\in X$. Then $y_{v'\otimes w', x'}=cy_{v\otimes w, x}$ for some $c\in\mathbb{C}^{\times}$ if and only if $v'\otimes w' - v\otimes w \in A^{*}$ and $x'=bx$ for some $b\in\mathbb{C}$.
\end{lemma}
\begin{proof}
If $y_{v'\otimes w',x'}=cy_{v\otimes w,x}$ for some $c\in\mathbb{C}^{\times}$, the the supports of the two functions are identical. Thus, $-(v\otimes w)+A^{*}=-(v'\otimes w')+A^{*}$ so that $v\otimes w - v'\otimes w' \in A^{*}$. Also, we have $y_{v'\otimes w', x'}(v'\otimes w')=cy_{v\otimes w, x}(v'\otimes w')$ which means
\begin{center}
$\chi(<<v'\otimes w',v'\otimes w'>>/2)x'=c\chi(<<v\otimes w, v'\otimes w'>>/2)x$

$x'=c\chi(<<v\otimes w, v'\otimes w'>>/2)x$
\end{center}

Conversely, assume $v'\otimes w' - v\otimes w \in A^{*}$ and $x'=bx$ for $b\in\mathbb{C}$. So $v'\otimes w'= v\otimes w +a^{*}$ for some $a^{*}\in A^{*}$. Then we have
\begin{align*}
y_{v'\otimes w', x'}(z)&=y_{v\otimes w+a^{*},bx}(z)\\
&=\chi(<<v\otimes w +a^{*}, z>>/2)bx\\
&=b\chi(<<a^{*},z>>/2)\chi(<<v\otimes w, z>>/2)x\\
&=b\chi(<<a^{*}, z>>/2)y_{v\otimes w, x}(z)
\end{align*}
\noindent for all $z\in -(v\otimes w)+A^{*}$. Hence we have $y_{v'\otimes w', x'}=cy_{v\otimes w, x}$ where $c=b\chi(<<a^{*}, z>>/2)$.  
\end{proof}

\begin{theorem}
With the notation as above, let $k$ be a positive integer and let $Y_{k}$ be the set of functions in $Y$ supported on $(M^{k})^{*}=P^{-k}A^{*}$. Then the following hold:

\begin{enumerate}
\item[(i)] $U^{2k+1}_{1}(\mathcal{A})$ and $U^{4k+2}_{2}(\mathcal{A})$ fix $Y_{k}$ pointwise

\item[(ii)] If $f$ is in $Y_{k}$ and $(h,1)\in U^{k}_{1}(\mathcal{A})$ (resp. $U^{2k}_{2}(\mathcal{A})$), then

\begin{center}
$\omega_{\chi}(h,1)f(v\otimes w)=\rho_{A}(2c(h)(v\otimes w))\chi(<<v\otimes w, c(h)(v\otimes w)>>)f(v\otimes w)$
\end{center}
\end{enumerate}
\end{theorem}
\begin{proof}
\begin{enumerate}
\item[(i)] Let $(g,1)\in U^{2k+1}_{1}(\mathcal{A})\cap SL_{2}(\mathcal{O}_{F})\rtimes \{1_{E}\}$. Then $(g,1)=g $ under the automorphism $\sigma$ and using Proposition $7$, Section $6$, we will show that $\rho_{A}(2c(g)(v\otimes w))$ and $\chi(<<v\otimes w, c(g)(v\otimes w)>>)$ are trivial. We take $v\otimes w\in(P^{-k}\Gamma)\otimes(\Gamma')^{*}$. We will start with the latter:
\begin{align*}
\chi(<<v\otimes w, c(g)(v\otimes w)>>)&=\chi(Tr(<v,c(g)v>\overline{<w,w>'}))\\
&=\chi(Tr(<v,c(g)v> N_{D/F}(w)))
\end{align*}
\noindent Now write $g=1+x$, $x\in P^{2k+1}$, i.e. $x=\mat{x_{1}\varpi^{2k+1}}{x_{2}\varpi^{2k+1}}{x_{3}\varpi^{2k+1}}{x_{4}\varpi^{2k+1}}$, then 
\begin{align*}
c(g)v&=(1-g)(1+g)^{-1}v\\
&=(-x)(2+x)^{-1}v\\
&=-2^{-1}x(1+2^{-1}x)^{-1}v\\
&=-2^{-1}x(1-2^{-1}x+(2^{-1}x)^{2}-(2^{-1}x)^{3}+...)v\\
&=-2^{-1}xv+(2^{-1}x)^{2}v-(2^{-1}x)^{3}v+...\\
&=\displaystyle{\sum^{\infty}_{i=1} (-1)^{i}(2^{-1}x)^{i}v}
\end{align*}
\noindent Hence the above trace formula will become 
\begin{center}
$\chi(Tr(<v,c(g)v>)N_{D/F}(w))=\chi(Tr(<v,\displaystyle{\sum^{\infty}_{i=1} (-1)^{i}(2^{-1}x)^{i}v}>N_{D/F}(w)))$
\end{center}
\noindent Note that the term with the smallest order in above expansion is $2^{-1}xv$ and hence $<v,-2^{-1}xv>N_{D/F}(w)$=$<a\varpi^{-k}{\bf{u}}+b\varpi^{-k}{\bf{v}},2^{-1}x(a\varpi^{-k}{\bf{u}}+b\varpi^{-k}{\bf{v}})>N_{D/F}(w)$ where ${\bf{u,v}}$ are the base vectors and $a,b\in\mathcal{O}_{E}$. 
\begin{align}
<v,-2^{-1}xv>&=<a\varpi^{-k}{\bf{u}}+b\varpi^{-k}{\bf{v}},2^{-1}x(a\varpi^{-k}{\bf{u}}+b\varpi^{-k}{\bf{v}})>\\
&=<a\varpi^{-k}{\bf{u}}+b\varpi^{-k}{\bf{v}},2^{-1}x(a\varpi^{-k}{\bf{u}}+b\varpi^{-k}{\bf{v}})>\\
&=2^{-1}(a\bar{a}\bar{x_{3}}\varpi+a\bar{x_{4}}\bar{b}\varpi-\bar{a}\bar{x_{1}}b\varpi-x_{2}b\bar{v}\varpi)
\end{align}

Then the valuation $v((a\bar{a}\bar{x_{3}}\varpi+a\bar{x_{4}}\bar{b}\varpi-\bar{a}\bar{x_{1}}b\varpi-x_{2}b\bar{v}\varpi)N_{D/F}(w))=1-1=0$ and since $\chi$ has a conductor $\mathcal{O}_{F}$, this term and all the following terms with bigger valuation will vanish.

To show that the representation $\rho_{A}(2c(g)(v\otimes w))=\rho_{A}((2c(g)v)\otimes w)$ is trivial, we consider again the valuation of involved terms: $v(c(g)v)=v(c(g))+v(v)\geq 0$ for $g\in U^{2k+1}_{1}(\mathcal{A})$ and hence representation will be trivial. A similar argument shows that $U^{4k+2}_{2}(\mathcal{A})$ fixes $Y_{k}$ pointwise. The different level is due to the ramification of $D$ over $F$, i.e. $g=1+x, x\in P^{4k+2}_{D} => x\in P^{2k+1}$. 

\vspace{1cm}
\item[(ii)] Again $(g,1)=g$ and we use Proposition $7$, Section $6$. $H_{M}=\{g\in G: (g-1)(M^{k})^{*}\subset A^{*}\}=\{(g-1)P^{-k}A^{*}\subset A^{*}\}$, it follows $g\in U^{k}_{1}(\mathcal{A})$ or $U^{2k}_{2}(\mathcal{A})$. 
\end{enumerate}
\end{proof}

\begin{theorem}
With the notation as in Theorem $9$.
\begin{enumerate}
\item[(i)] Let $(h,1)\in U^{k+1}_{1}(\mathcal{A})$. If $f$ is in $Y_{k}$, then $f$ transforms according to $\psi_{b_{1}}$ under the actions of $U^{k+1}_{1}(\mathcal{A})$ where $b_{1}\in E^{0}$. 
\begin{center}
$b_{1}=-\frac{\varpi^{-k+2}}{2}N_{D/F}(w)\mat{-\varpi^{-k-1}(\bar{a}b+a\bar{b})}{2\varpi^{-k-1}a\bar{a}}{-2\varpi^{-k-1}b\bar{b}}{\varpi^{-k-1}(a\bar{b}+\bar{a}b)}$
\end{center}

\item[(ii)] Let $(h,1)\in U^{2k+2}_{2}(\mathcal{A})$. If $f$ is in $Y_{k}$, then $f$ transforms according to $\psi_{b_{2}}$ under the actions of $U^{2k+2}_{2}(\mathcal{A})$ where $b_{2}\in D^{0}$.
\begin{center}
$b_{2}=-\frac{\varpi}{2}N_{D/F}(w)(a\bar{b}-\bar{a}b)$ 
\end{center}
\end{enumerate}
\end{theorem}
\begin{proof}
\begin{enumerate}
\item[(i)] 
Let $(h,1)=h\in U^{k+1}_{1}(\mathcal{A}), h=1+x, x\in P^{k+1}$ and take $v\otimes w\in (P^{-k}\Gamma)\otimes (\Gamma')^{*}$. Consider the valuation of involved terms in $\rho_{A}((2c(h)v)\otimes w)$. As in the Theorem $9$(i), the term with the smallest order in an expansion is $2^{-1}xv$ and hence its valuation $v(xv)=v(x)+v(v)=k+1-k=1 \geq 0$. It follows that the $\rho_{A}$ action is trivial and 
\begin{center}
$\omega_{\chi}(h,1)f(v\otimes w)=\chi(<<v\otimes w, (c(h)v)\otimes w>>)f(v\otimes w)$
\end{center}
\noindent Now arguing as in Theorem $9$ (i), we get 
\begin{align*}
\chi(<<v\otimes w, (c(h)v)\otimes w>>)
&=\chi(Tr(<v,c(h)v>N_{D/F}(w)))\\
&=\chi(Tr(<v,\displaystyle{\sum_{i=1}^{\infty}} (-1)^{i}(2^{-1}x)^{i}v> N_{D/F}(w)))
\end{align*}
\noindent which is the same as 
\begin{center}
$\chi(Tr(<v,-2^{-1}xv>N_{D/F}(w)+<v,\displaystyle{\sum_{i=2}^{\infty}} (-1)^{i}(2^{-1}x)^{i}v>N_{D/F}(w)))$.
\end{center}
\noindent The second term in above is in $P_{E}$ and since $\chi$ has a conductor $\mathcal{O}_{F}$, we obtain
\begin{center}
$\chi(Tr(<v,\displaystyle{\sum_{i=2}^{\infty}} (-1)^{i}(2^{-1}x)^{i}v>N_{D/F}(w)))=1$
\end{center}
\noindent and hence
\begin{align}
\chi(<<v\otimes w, (c(h)v)\otimes w>>)=\chi(Tr(<v,-2^{-1}xv>N_{D/F}(w)))
\end{align}

\noindent We will look at the computations separately. First, we want to explicitly compute $<v, -2^{-1}xv>$. 
\begin{center}
$xv= \mat{x_{1}\varpi^{k+1}}{x_{2}\varpi^{k+1}}{x_{3}\varpi^{k+1}}{x_{4}\varpi^{k+1}}\matt{\varpi^{-k}a}{\varpi^{-k}b}$
\end{center}

\noindent which is precisely $(\varpi ax_{1}+\varpi bx_{2}){\bf{u}}+(\varpi ax_{3}+\varpi bx_{4}){\bf{v}}$. \\

\noindent Using this in above quadratic form will give us the explicit expression for $<v, -2^{-1}xv>$:
\begin{align}
<v, -2^{-1}xv>
&= <\varpi^{-k}a{\bf{u}}+\varpi^{-k}b{\bf{v}}, -2^{-1}(\varpi ax_{1}+\varpi bx_{2}){\bf{u}}+(\varpi ax_{3}+\varpi bx_{4}){\bf{v}}>\\
&=-2^{-1}\varpi^{-k+1}(x_{3}a\bar{a}+x_{4}a\bar{b}-x_{1}\bar{a}b-x_{2}b\bar{b})
\end{align}

\noindent Using the result from equation $(6)$ in equation $(4)$, we obtain the following:
\begin{align}
\chi(<<v\otimes w, (c(h)v)\otimes w>>)&=
\end{align}
\begin{align*}
&=\chi(Tr(-2^{-1}\varpi^{-k+1}(x_{3}a\bar{a}+x_{4}a\bar{b}-x_{1}\bar{a}b-x_{2}b\bar{b})N_{D/F}(w)))\\
&=\chi(-2^{-1}\varpi^{-k+1}N_{D/F}(w)(2x_{3}a\bar{a}+x_{4}a\bar{b}+x_{4}\bar{a}b-x_{1}\bar{a}b-x_{1}a\bar{b}-2x_{2}b\bar{b})).
\end{align*}

\noindent This formula corresponds to the trace of the following element:
\begin{center}
$-2^{-1}\varpi^{-k+1}N_{D/F}(w)\scriptsize{\mat{-\varpi^{-k-1}(\bar{a}b+a\bar{b})}{2\varpi^{-k-1}a\bar{a}}{-2\varpi^{-k-1}b\bar{b}}{\varpi^{-k-1}(a\bar{b}+\bar{a}b)}\mat{x_{1}\varpi^{k+1}}{x_{2}\varpi^{k+1}}{x_{3}\varpi^{k+1}}{x_{4}\varpi^{k+1}}}$
\end{center}

\noindent Now putting the equation $(4)$ and $(7)$ together, we are able to explicitly write the formula for the traceless element $b_{1}$ involved the theta correspondence:
\begin{align}
\chi(-2^{-1}\varpi^{-k+1}N_{D/F}(w)(2x_{3}a\bar{a}+x_{4}a\bar{b}+x_{4}\bar{a}b-x_{1}\bar{a}b-x_{1}a\bar{b}-2x_{2}b\bar{b}))=
\end{align}
\begin{align*}
&=\chi(Tr(-2^{-1}\varpi^{-k+1}N_{D/F}(w)\scriptsize{\mat{-\varpi^{-k-1}(\bar{a}b+a\bar{b})}{2\varpi^{-k-1}a\bar{a}}{-2\varpi^{-k-1}b\bar{b}}{\varpi^{-k-1}(a\bar{b}+\bar{a}b)}\mat{x_{1}\varpi^{k+1}}{x_{2}\varpi^{k+1}}{x_{3}\varpi^{k+1}}{x_{4}\varpi^{k+1}}}))\\
&=\psi(Tr(-\frac{\varpi^{-k+2}}{2}N_{D/F}(w)\mat{-\varpi^{-k-1}(\bar{a}b+a\bar{b})}{2\varpi^{-k-1}a\bar{a}}{-2\varpi^{-k-1}b\bar{b}}{\varpi^{-k-1}(a\bar{b}+\bar{a}b)} x))\\
&=\psi(Tr(b_{1}(x)))\\
&=\psi(Tr(b_{1}(h-1)))\\
&=\psi_{b_{1}}(h)
\end{align*}

\noindent where $b_{1}=-\frac{\varpi^{-k+2}}{2}N_{D/F}(w)\mat{-\varpi^{-k-1}(\bar{a}b+a\bar{b})}{2\varpi^{-k-1}a\bar{a}}{-2\varpi^{-k-1}b\bar{b}}{\varpi^{-k-1}(a\bar{b}+\bar{a}b)}$. 

\noindent It is clear that the element $b_{1}$ is traceless. 

\smallskip
\item[(ii)] With the notation as above, we take $v\otimes w\in (\Gamma)\otimes (P^{-k}(\Gamma')^{*})$, i.e. $v=a{\bf{u}}+b{\bf{v}}, w=\varpi^{-k}c+\varpi^{-k-1}d\delta$ where $a,b,c,d\in\mathcal{O}_{E}$. Now let $(h,1)\in U^{2k+2}_{2}(\mathcal{A})$,i.e. $(h,1)=h=1+x, x\in P^{2k+2}_{D}$. 

Consider the valuation of involved terms in $\rho_{A}(v\otimes (2c(h)w))$. As in Theorem $9$(i), the term with the smallest order in an expansion is $2^{-1}xw$ and hence its valuation $v(xw)=v(x)+v(w)=k+1-k-1=0$. It follows that $\rho_{A}$ is trivial and 
\begin{center}
$\omega_{\chi}(h,1)f(v\otimes w)=\chi(<<v\otimes w, v\otimes (c(h)w)>>)f(v\otimes w)$
\end{center} 

Hence we have
\begin{align}
\chi(<<v\otimes w, v\otimes (c(h)w)>>)&=\chi(Tr(<v,v>\overline{<w,c(h)w>'}))\\
&=\chi(Tr((a\bar{b}-\bar{a}b)\overline{1/2 Tr_{D/E}(w\overline{c(h)w})}))
\end{align}

Since $\chi$ has conductor $\mathcal{O}_{F}$, all terms vanish in $Tr_{D/E}$ will vanish but the first one. Therefore, we obtain

\begin{align}
\chi(Tr((a\bar{b}-\bar{a}b)\overline{1/2 Tr_{D/E}(w(-2^{-1})\bar{x}\bar{w})}))=
\end{align}

\begin{align*}
&=\chi(Tr((a\bar{b}-\bar{a}b)(-1/4\overline{N_{D/F}(w)Tr_{D/E}(\bar{x})})))\\
&=\chi(Tr((a\bar{b}-\bar{a}b)(-1/2 N_{D/F}(w)(x))))\\
&=\psi(Tr(-\frac{\varpi}{2}(a\bar{b}-\bar{a}b)N_{D/F}(w)(h-1)))\\
&=\psi_{b_{2}}(h)
\end{align*}

\noindent where $b_{2}=-\frac{\varpi}{2}(a\bar{b}-\bar{a}b)N_{D/F}(w)$, clearly $b_{2}$ is traceless and an element of $D$.
\end{enumerate}
\end{proof} 

\begin{theorem}
With the notation as above, $N(b_{2})=\det b_{1}$.
\end{theorem}
\begin{proof}
In Theorem $1$, $w\in (\Gamma')^{*}, w=c+\varpi^{-1}d\delta, c,d\in\mathcal{O}_{E}$ and hence $N_{D/F}(w)=N_{D/F}(c)-\varpi^{-2}N_{D/F}(d)\delta^{2}$. In Theorem $2$, $w'\in (P^{-k}(\Gamma')^{*}), w'=\varpi^{-k}c+\varpi^{-k-1}d\delta$ and hence $N_{D/F}(w')=\varpi^{-2k}N_{D/F}(c)-\varpi^{-2k-2}N_{D/F}(d)\delta^{2}=\varpi^{-2k}N_{D/F}(w)$. So in fact, these two norms differ by a term $\varpi^{-2k}$ which we will factor out in computations for $N(b_{2})$. 
\begin{align*}
\det b_{1}&=\det (-\frac{\varpi^{-k+2}}{2}N_{D/F}(w)\mat{-\varpi^{-k-1}(\bar{a}b+a\bar{b})}{2\varpi^{-k-1}a\bar{a}}{-2\varpi^{-k-1}b\bar{b}}{\varpi^{-k-1}(a\bar{b}+\bar{a}b)})\\
&=-\frac{\varpi^{-4k+2}}{4}N_{D/F}^{2}(w)(\bar{a}b+a\bar{b})^{2}+\varpi^{-4k+2}N_{D/F}^{2}(w)a\bar{a}b\bar{b}\\
&=-\frac{\varpi^{-4k+2}}{4}N_{D/F}^{2}(w)((a\bar{b})^{2}+2a\bar{a}b\bar{b}+(\bar{a}b)^{2}-4a\bar{a}b\bar{b})\\
&=-\frac{\varpi^{-4k+2}}{4}N_{D/F}^{2}(w)(a\bar{b}-\bar{a}b)^{2}
\end{align*}
The right hand side is equal to
\begin{align*}
N(b_{2})&=N(-\frac{\varpi}{2}(a\bar{b}-\bar{a}b)N_{D/F}(w'))\\
&=N(-\frac{\varpi}{2}(a\bar{b}-\bar{a}b)\varpi^{-2k}N_{D/F}(w))\\
&=(-\frac{\varpi^{-2k+1}}{2}(a\bar{b}-\bar{a}b)N_{D/F}(w))\overline{(-\frac{\varpi^{-2k+1}}{2}(a\bar{b}-\bar{a}b)N_{D/F}(w))}\\
&=-\frac{\varpi^{-4k+2}}{4}N_{D/F}^{2}(w)(a\bar{b}-\bar{a}b)^{2}
\end{align*}
\end{proof}

\begin{corollary}
Thus with the notation as above, $b_{1}$ and $b_{2}$ belong to corresponding conjugacy classes in $GL_{2}(F)$ and in $D^{\times}$.
\end{corollary}
For the correspondence between conjugacy classes in $GL(n)$ and division algebra, see\cite{Rogawski}.

\end{document}